\documentclass[final]{siamart171218}

\newsiamremark{notation}{Notation}
\newsiamremark{remark}{Remark}

\headers{Global supersolutions for systems with taxis}{A. Zhigun}

\title{Generalised global supersolutions with  mass control for systems with taxis}

\author{Anna Zhigun\thanks{Queen's University Belfast, School of Mathematics and Physics,
University Road, Belfast BT7 1NN, Northern Ireland, UK 
  (\email{A.Zhigun@qub.ac.uk}).}}

\usepackage{amsopn}

\ifpdf
\hypersetup{
  pdftitle={Generalised global supersolutions with  mass control for systems with taxis},
  pdfauthor={A. Zhigun}
}
\fi

\usepackage{amsfonts,amsmath}
\usepackage{empheq}

\usepackage{cases}
\usepackage{mathabx}
\usepackage{bbm}
\usepackage{xfrac}
\usepackage{fancyhdr}
\usepackage{constants}
\newcommand{\parenthezises}[1]{\arabic{#1}}
\newconstantfamily{C}{
symbol=C,
format=\parenthezises,
reset={section}
}
\newconstantfamily{M}{
symbol=M,
format=\parenthezises,
reset={section}
}
\newconstantfamily{B}{
symbol=B,
format=\parenthezises,
reset={section}
}
\usepackage{enumerate}
\usepackage{graphicx}
\graphicspath{ {images/} }
\usepackage{wrapfig}
\usepackage{figbib}

\usepackage{mathtools}
\allowdisplaybreaks
\begin{document}

\newcommand{\D}{\mathbb{D}}
\newcommand{\R}{\mathbb{R}}
\newcommand{\N}{\mathbb{N}}

\def\diam{\operatorname{diam}}
\def\dist{\operatorname{dist}}
\def\diver{\operatorname{div}}
\def\ess{\operatorname{ess}}
\def\inner{\operatorname{int}}
\def\osc{\operatorname{osc}}
\def\sign{\operatorname{sign}}
\def\supp{\operatorname{supp}}
\newcommand{\BMO}{BMO(\Omega)}
\newcommand{\LOne}{L^{1}(\Omega)}
\newcommand{\LOnen}{(L^{1}(\Omega))^d}
\newcommand{\LTwo}{L^{2}(\Omega)}
\newcommand{\Lq}{L^{q}(\Omega)}
\newcommand{\Lp}{L^{2}(\Omega)}
\newcommand{\Lpn}{(L^{2}(\Omega))^d}
\newcommand{\LInf}{L^{\infty}(\Omega)}
\newcommand{\HOneO}{H^{1,0}(\Omega)}
\newcommand{\HTwoO}{H^{2,0}(\Omega)}
\newcommand{\HOne}{H^{1}(\Omega)}
\newcommand{\HTwo}{H^{2}(\Omega)}
\newcommand{\HmOne}{H^{-1}(\Omega)}
\newcommand{\HmTwo}{H^{-2}(\Omega)}

\newcommand{\LlogL}{L\log L(\Omega)}

\def\avint{\mathop{\,\rlap{-}\!\!\int}\nolimits}


\maketitle
\begin{abstract}
 The existence of  generalised global supersolutions with a control upon the total muss is established for a wide family of  parabolic-parabolic chemotaxis systems and general integrable initial data in any space dimension. It is verified that as long as a  supersolution  of this sort remains smooth, it coincides with the  classical solution. At the same time, the proposed construction provides solvability beyond a blow-up time.  
 The considered class of systems includes the basic form of the Keller-Segel model as well as the case of a chemorepellent.
\end{abstract}

\section{Introduction}
Ever since it was introduced \cite{KS70,KS71} in the 1970's, the Keller-Segel model for chemotaxis and its modifications have been extensively studied, compare, e.g., reviews \cite{Horstmann,BBTW} on the available results. The  best-known model representative has the form
\begin{subequations}\label{KScl}
 \begin{alignat}{3}
  &\partial_tu=\nabla\cdot\left(\nabla u-\chi u\nabla v\right)\quad&&\text{ in }\R^+\times\Omega,\\
  &\partial_t v=\Delta v-v+u\quad&&\text{ in }\R^+\times\Omega,\\
  &\partial_{\nu} u=\partial_{\nu} v=0\quad&&\text{ in }\R^+\times\partial\Omega,\\
  &u(0,\cdot)=u_0,\ v(0,\cdot)=v_0\quad&&\text{ in }\Omega,
 \end{alignat}
\end{subequations}
where $\Omega$ is a smooth bounded domain in $\R^N$, $N\in\N$, with the corresponding outer normal unit vector $\nu$ on $\partial\Omega$. The constant $\chi>0$ is the so-called chemotactic sensitivity coefficient.

It has long been suspected that in higher dimensions  \cref{KScl} possesses solutions which blow-up in finite time. The later means that there are solutions which start as classical at $t=0$ but exist only up to some finite time  $T_{max}\in\R^+$ when the aggregation and collapse of the $u$-component occur, i.e., \begin{align*}
\underset{t\uparrow T_{max}}{\lim\sup}\, \|u(t,\cdot)\|_{L^{\infty}(\Omega)}=\infty.                                                                                                                                                                                                                                                                                                                                                                  \end{align*}
This hypothesis  was finally confirmed in  {\cite{Winkler2013}} for the case when $\Omega$ is a ball. 

One of the  questions naturally arising in connection with  the blow-up phenomena is the following: how can a solution be reasonably extended  beyond its blow-up time? It is clear that the only nonlinear term in the system, $\chi u\nabla v$, creates difficulties. Indeed, it is not obvious as to how the product between $u$ and $\nabla v$ should be interpreted if, for instance, already 
\begin{align*}
 \underset{t\uparrow T_{max}}{\lim\sup}\, \|uv(t,\cdot)\|_{L^1(\Omega)}=\infty
\end{align*}
can occur in the blow-up case {\cite{Winkler2013}}. 
This issue has been addressed with success for some  parabolic-elliptic modifications of \cref{KScl}  in dimension $N=2$. Usually one introduces  a suitable family of regular problems parametrised by a small parameter $\varepsilon>0$ in such a way that by taking $\varepsilon=0$  the original system is formally recovered. 
The main challenges are then to find the correct limit problem for $\varepsilon\rightarrow 0$ and  to rigorously prove convergence. Formal asymptotics undertaken in   \cite{Velazquez2004_1,Velazquez2004_2}  for one class of such regularisations revealed the structure of the blow-up solutions and indicated what the corresponding limit problem should be. The local in time well-posedness of the later was also verified \cite{Velazquez2004_3}. 
Several subsequent works  deal with the construction of  generalised weak solutions in $\R^2$ \cite{DolbSchm2009} or a bounded domain \cite{LuckSugVel} and establish convergence to it. An important mark of the resulting limit problems  is that they involve measure-valued terms, which, moreover, turn out to be dependent upon the choice of a specific regularisation. Measures necessarily appear in those parts of the solution formulation where integrability can be lost due to a potential blow-up. The  solution concepts in both \cite{DolbSchm2009} and \cite{LuckSugVel}  rely heavily upon the properties of the Green function of the Laplace operator in planar domains, as well as  a certain symmetry of the model: both the coefficient before $\nabla v$ in the convective flux in the equation for $u$ and the production term in the $v$-equation are linear in $u$. It has remained open as to how these constructions could be extended to the cases in which the space dimensions $N$ is larger than two, and/or the system is fully parabolic,  and/or $\chi$ depends upon $u$. The present study aims to fill this gap.

System \cref{KScl} can also be brought into  comparison with its modification involving a  logarithmic chemotactic sensitivity function: in this case the term $\chi u\nabla v$ is replaced by $\chi u\nabla \ln(v)=\chi \frac{u}{v}\nabla v$. Owing to the damping effect of $\frac{1}{v}$ 
at large values of $v$, the corresponding solutions are less prone to formation of strong singularities. Very recently a new, generalised, solution concept has been introduced \cite{LankWink2017} and further developed \cite{ZhigunKSLog} which makes an extensive use of this property. { The} underlying idea { was} to replace the equation for $u$ by a variational inequality which describes the evolution of a certain coupling of both $u$ and $v$. There results a supersolution concept with a mass control from above, and it has been shown that such a supersolution, should it be regular, actually coincides with the  classical solution. By switching to an inequality, one avoids the necessity to include a non-negative measure-valued term, of which, other than in the cases studied in \cite{DolbSchm2009,LuckSugVel}, not much is presently known. As to the employed  coupled quantities, they enjoy better regularity than $u$ generally does. It ultimately allows to  sustain some sort of  solvability even if a blow-up does take place in some parts of the system. This approach has previously been used in the context of degenerate haptotaxis systems \cite{ZSH,ZSU}. 

In this paper we elaborate a global solvability concept { for  the following generalisation of \cref{KScl}:}
\begin{subequations}\label{KS}
 \begin{alignat}{3}
  &\partial_tu=\nabla\cdot\left(\nabla u-c(u,v)\nabla v\right)\quad&&\text{ in }\R^+\times\Omega,\label{Equ}\\
  &\partial_t v=\Delta v-v+u\quad&&\text{ in }\R^+\times\Omega,\label{Eqv}\\
  &\partial_{\nu} u=\partial_{\nu} v=0\quad&&\text{ in }\R^+\times\partial\Omega,\label{bc}\\
  &u(0,\cdot)=u_0,\ v(0,\cdot)=v_0\quad&&\text{ in }\Omega,
 \end{alignat}
\end{subequations}
where 
\begin{align*}
 &c(u,v):=\chi(u,v)u
\end{align*}
for an arbitrary  
\begin{align}
 {\chi\in C(\R_0^+\times \R_0^+).}\nonumber
\end{align}
It is important to note that no further restrictions upon $\chi$ are required. In fact, just as it is the case for the solvability in the classical sense, the upper-triangular structure of system \cref{KS} turns out to be a decisive factor for the generalised solvability as well. More precisely, the outcome of this work can be summarised as follows: the classical theory  for upper-triangular systems (see, e.g., \cite{Amann1}) implies that \cref{KS} possesses a { solution which is classical} as long as it doesn't blow-up, and if it does, then it continuous to exist as generalised in terms of {\cref{defsol}} given below.  

Our construction is based on the ideas developed in   \cite{LankWink2017,ZhigunKSLog} for the case of a logarithmic sensitivity (which is not covered by this work since we assume $\chi$ to be continuous). The main difference lies in the choice of suitable couplings of $u$ and $v$: we use cut-offs in order to overcome possible unboundedness of any of the two components. 
As in \cite{LankWink2017,ZhigunKSLog}, the  generalised solution is obtained by means of a limit procedure for an   approximation family which { in the present case} involves a regularisation  
{ in the chemotactic coefficient} (see \cref{KSe} below). However, other { regularisations are} also possible. It is to expect that just as in \cite{DolbSchm2009,LuckSugVel} the resulting solution depends upon the choice of a particular approximation.  

The rest of the paper is organised in the following way. We begin with a preliminary {\cref{SecPrelim}} where we fix some notations which will be used throughout and also recall some results involving measures. 
In {\cref{SecSol}} we introduce our supersolution concept for \cref{KS} (cl. {\cref{defsol}}), and   two main results: {\cref{TheoEx}} deals with existence of such supersolutions, while {\cref{classol}} establishes a link to the classical solvability. In {\cref{secreg}} we prepare some  ingredients which are necessary to our proof of  {\cref{TheoEx}} in {\cref{SecEx}}. {\cref{classol}} is proved in the closing {\cref{SecCl}}. 

\subsection*{Acknowledgement}
The author expresses her gratitude to the reviewers for their helpful comments.
\section{Preliminaries}\label{SecPrelim}
A number of notational conventions are used in this paper for the sake of conciseness:
\begin{notation}\label{Notation}~
\begin{enumerate}[(1)]
 \item For any index $i$, a quantity $C_i$ denotes a positive constant or function.
 \item 
Dependence upon such parameters as: the space { dimension $N$}, domain $\Omega$, function $c$,  the $L^1$-norms of the initial conditions $u_0$ and $v_0$, as well as constants $a$ and $b$ (are introduced below) is mostly 
{\bf not} indicated in an explicit way.
\item\label{cut} { Given $n,l\in\N$, we introduce the cut-offs}
 \begin{align}
 &\overline{u}:=(n-u)_+,\qquad \overline{v}:=(l-v)_+. \nonumber
\end{align}
Terms $\overline{u_k}, \overline{u_{k_m}},\overline{v_k}$, $\overline{v_{k_m}}$  are to be understood in the same fashion. 
\end{enumerate}
\end{notation}
\subsection{Functional spaces}\label{SecFSp}
We assume the reader to be familiar with the standard spaces of continuous and continuously differentiable functions, as well as the Lebesgue, Sobolev,    and Bochner spaces and standard results concerning them. 

As usual, ${\cal M}(\overline{\Omega})$ denotes the Banach space of signed Radon measures in $\Omega$, while   ${\cal M}_+(\overline{\Omega})$ stands for the set of all positive Radon measures. 

We also make use of the Banach space
\begin{align*}
L^{\infty}_{w-*}(\R^+;{\cal M}(\overline{\Omega})):=&\left\{\mu:\R^+\rightarrow {\cal M}(\overline{\Omega})\quad\text{is weak}-*\text{ measurable and}\right.\\ 
&\ \left. \left\|\mu\right\|_{L^{\infty}_{w-*}(\R^+;{\cal M}(\overline{\Omega}))}:=\left\|\|\mu(\cdot)\|_{{\cal M}(\overline{\Omega})}\right\|_{L^{\infty}(\R^+)}<\infty\right\}
\end{align*}
and its closed subset
\begin{align*}
&L^{\infty}_{w-*}(\R^+;{\cal M}_+(\overline{\Omega})):=\left\{\mu:\R^+\rightarrow {\cal M}_+(\overline{\Omega})|\ \mu\in L^{\infty}_{w-*}(\R^+;{\cal M}(\overline{\Omega}))\right\}.
\end{align*}
 We identify functions which coincide a.e. in $\R^+$. It is known \cite[sections 8.18.1-8.18.2]{Edwards} that
$L^{\infty}_{w-*}(\R^+;{\cal M}(\overline{\Omega}))$ is isometrically isomorphic to the continuous dual of the Bochner space $L^1(\R^+{;}C(\overline{\Omega}))$ via the duality paring 
\begin{align*}
 \left<\mu,\varphi\right>=\int_{\R^+}\int_{\overline{\Omega}}\varphi(t)(x)\, (d\mu(t))(x)dt.
\end{align*}
Since $L^1(\R^+{;}C(\overline{\Omega}))$ is separable, the Banach-Alaoglu theorem implies that { closed} balls in $L^{\infty}_{w-*}(\R^+;{\cal M}(\overline{\Omega}))$ are weak-$*$ sequentially compact.

\subsection{Divergence-measure fields and their normal traces}\label{AppA}
In this sequel some facts concerning the  divergence-measure fields are  { presented in the same way as  in \cite[Appendix A]{ZhigunKSLog}}. We use them in order to give weak reformulations of the Neumann boundary conditions (see {\cref{defsol}} below). 


  Recall the definition of the Banach space of divergence-measure fields \cite{ChenFrid2001} and its norm: 
\begin{align}
&{\cal DM}^p(\Omega):=\{F\in (L^p(\Omega))^N|\ \ \nabla\cdot F\in {\cal M}(\overline{\Omega})\},\nonumber\\
 &\|F\|_{{\cal DM}^p(\Omega)}:=\|F\|_{(L^p(\Omega))^N}+\|\nabla\cdot F\|_{{\cal M}(\overline{\Omega})}.\nonumber
\end{align}
Here we assume that $p\in\left(1,\frac{N}{N-1}\right)$, which is sufficient for our needs. 
Following \cite{ChenFrid2001} we introduce a generalisation of the normal trace over the boundary of $\partial\Omega$ which automatically satisfies a  Gauss-Green formula:  
\begin{align}
 \left<F\cdot\nu|_{\partial\Omega},\varphi\right>:=\int_{\Omega}F\cdot\nabla ({\cal E}\varphi)\,dx+\int_{\overline{\Omega}} ({\cal E}\varphi)\, d(\nabla\cdot F)\qquad\text{for all }\varphi\in  W^{\frac{1}{p},\frac{p}{p-1}}(\partial\Omega).\label{GG}
\end{align}
Here ${\cal E}: W^{\frac{1}{p},\frac{p}{p-1}}(\partial\Omega)\rightarrow W^{1,\frac{p}{p-1}}(\Omega)$ is a usual   extension operator, i.e., a continuous right inverse of the  corresponding trace operator. 
  It is known (see \cite[Theorem 2.1]{ChenFrid2001}) that  $F\cdot\nu|_{\partial\Omega}\in W^{-\frac{1}{p},p}(\partial\Omega)$ and  doesn't depend upon the particular choice of ${\cal E}$. Formula \cref{GG} ensures the following implication:
  \begin{align}
   \begin{rcases}F_n\underset{m\rightarrow\infty}{{\rightharpoonup}}F&\text{ in }(L^p(\Omega))^N,\\ \nabla \cdot F_n\underset{m\rightarrow\infty}{\overset{*}{\rightharpoonup}}\nabla\cdot F&\text{ in }{\cal M}(\overline{\Omega})
   \end{rcases}\quad\Rightarrow\quad  F_n\cdot\nu|_{\partial\Omega}\underset{m\rightarrow\infty}{\overset{*}{\rightharpoonup}}F\cdot\nu|_{\partial\Omega}\quad\text{in }W^{-\frac{1}{p},p}(\partial\Omega).\label{contitr}
  \end{align}


\section{Generalised supersolutions to \texorpdfstring{\cref{KS}}{}}\label{SecSol}
We make the following general assumptions on the data involved:
\begin{align}
&{ c(u,v)=\chi(u,v)u,\qquad u,v\in \R_0^+,}\label{cchi}\\
&\chi\in C(\R_0^+\times \R_0^+),\label{chiC}\\
 &u_0,v_0\in L^1(\Omega),\label{uv0L1}\\
 &0\leq u_0,v_0\not\equiv0.\label{uvp}
\end{align}
Motivated by an idea from  \cite{LankWink2017} and our previous work \cite{ZhigunKSLog}, we introduce a  concept of a generalised supersolution with a mass control. 
\begin{definition}[Generalised supersolution]\label{defsol} { Let $c$ and $(u_0,v_0)$ satisfy  \cref{cchi}-\cref{uvp}}.  We call a pair of measurable functions $(u,v):\R_0^+\times\overline{\Omega}\rightarrow {\R_0^+}\times {\R^+}$ a \underline{generalised} \underline{supersolution} to system \cref{KS} if there exist some:
\begin{enumerate}[(1)]
\item element $\mu\in L^{\infty}_{w-*}(\R^+;{\cal M}_+(\overline{\Omega}))$,
 \item numbers $a,b>2$,
 \item function $M:\N\times\N\rightarrow \R^+_0$, { $M=M(n,l)$,}
\end{enumerate}
such that { for all $n,l\in\N$}:
\begin{enumerate}[(i)]
  \item\label{defi} $u\in L^{\infty}(\R^+;L^1(\Omega))$,\  $v\in L^1_{loc}(\R_0^+;W^{1,1}(\Omega))$,\ $v^{-1}\in L_{loc}^{\infty}(\R^+\times\overline{\Omega})$;
  \item $\overline{u}^{a}\overline{v}^{b},\ \overline{v}^{b}\in L^2_{loc}(\R^+_0;H^1(\Omega))$;
  \item \label{inDMp} $\int_0^{\infty}\psi\left(\left(\nabla\left( \overline{u}^{a}\overline{v}^{b}\right)-\overline{u}^{a}\nabla \overline{v}^{b}\right)+\frac{1}{b}\left(b\left({ M} +\overline{u}^{a}\right)-ac(u,v)\overline{u}^{a-1}\overline{v}\right)\nabla \overline{v}^{b}\right)\,ds$,\\   $\int_0^{\infty}\psi \nabla v\,ds\in {\cal DM}^{p}(\Omega)$   for all $\psi\in C^1_0(\R_0^+)$ for some $p\in\left(1,\frac{N}{N-1}\right)$;
\item for all $0\leq\varphi\in C^1(\overline{\Omega})$ and $0\leq\psi\in C^1_0(\R_0^+)$ it holds that 
\begin{align}
&-\int_0^{\infty}\partial_t\psi\int_{\Omega}\left({ M} +\overline{u}^{a}\right)\overline{v}^{b}\varphi\,dxdt-\psi(0)\int_{\Omega}\left({ M} +\overline{u_0}^{a}\right)\overline{v_0}^{b}\,dx\nonumber\\
 \leq&\int_0^{\infty}\psi\int_{\Omega}-4\left(\frac{a-1}{a}\left|\nabla \left(\overline{u}^{\frac{a}{2}}\overline{v}^{\frac{b}{2}}\right)-\overline{u}^{\frac{a}{2}}\nabla \overline{v}^{\frac{b}{2}}\right|^2\right.\nonumber\\
 &\phantom{aaaaaa}\left.+\frac{1}{b}\overline{u}^{\frac{a}{2}-1}\left(2b\overline{u}-(a-1)c(u,v)\overline{v}\right)\left(\nabla \left(\overline{u}^{\frac{a}{2}}\overline{v}^{\frac{b}{2}}\right)-\overline{u}^{\frac{a}{2}}\nabla \overline{v}^{\frac{b}{2}}\right)\cdot\nabla \overline{v}^{\frac{b}{2}}\right.\nonumber\\
 &\phantom{aaaaaa}\left.+\frac{1}{b}\left((b-1)\left({ M}+\overline{u}^{a}\right)-ac(u,v)\overline{u}^{a-1}\overline{v}\right)\left|\nabla \overline{v}^{\frac{b}{2}}\right|^2\right)\varphi\nonumber\\
 &\phantom{aaaaaa}-\left(\left(\nabla\left( \overline{u}^{a}\overline{v}^{b}\right)-\overline{u}^{a}\nabla \overline{v}^{b}\right)+\frac{1}{b}\left(b\left({ M} +\overline{u}^{a}\right)-ac(u,v)\overline{u}^{a-1}\overline{v}\right)\nabla \overline{v}^{b}\right)\cdot\nabla\varphi\nonumber\\
 &\phantom{aaaaaa}-b\left({ M} +\overline{u}^{a}\right)\overline{v}^{b-1}\left(-v+u\right)\varphi\,dxdt;\label{superu}
 \end{align}
\item for all { $\varphi\in C^1(\overline{\Omega})$ and $\psi\in C^1_0(\R_0^+)$} it holds that
\begin{align}
 -\int_0^{\infty}\partial_t\psi\int_{\Omega}v\varphi\,dxds-\psi(0)\int_{\Omega}v_0\varphi\,dx=&\int_0^{\infty}\psi\int_{\Omega}-\nabla v\cdot\nabla\varphi+(-v+u)\varphi\,dx\nonumber\\
 &+\int_{\overline{\Omega}}\varphi \,(d\mu(s))(x)ds,\label{superv}
\end{align}
and 
\begin{align}
 \|u(t,\cdot)\|_{L^1(\Omega)}+\int_{\overline{\Omega}} \,(d\mu(t))(x)= \|u_{0}\|_{L^1(\Omega)}\qquad\text{for a.a. }t>0;\label{massu0}
\end{align}
\item
for all $\psi\in C^1_0(\R_0^+)$ it holds in $W^{-\frac{1}{p},p}(\partial\Omega)$ that 
\begin{align}
&\int_0^{\infty}\psi\left(\left(\nabla\left( \overline{u}^{a}\overline{v}^{b}\right)-\overline{u}^{a}\nabla \overline{v}^{b}\right)+\frac{1}{b}\left(b\left({ M} +\overline{u}^{a}\right)-ac(u,v)\overline{u}^{a-1}\overline{v}\right)\nabla \overline{v}^{b}\right)\,ds\cdot\nu|_{\partial\Omega}\nonumber\\
=&0,\label{prweakbc}\\
 &\int_0^{\infty}\psi\nabla v\,ds\cdot\nu|_{\partial\Omega}=0.\label{vweakbc}
\end{align}
\end{enumerate}
 Here $\bar u$ and $\bar v$ are as defined in {\cref{Notation}\eqref{cut}}, i.e., they implicitly depend upon $n$ and $l$. 
\end{definition}
\begin{remark}[Boundary conditions] The variational reformulations \cref{prweakbc}-\cref{vweakbc} of the boundary conditions \cref{bc} are consistent with the regularity assumptions in \eqref{inDMp}, cl. \cite[Theorem 2.1]{ChenFrid2001}. 
\end{remark}

Our result on existence now reads:
\begin{theorem}[Existence of generalised supersolutions]\label{TheoEx}
   Let $c$ and $(u_0,v_0)$ satisfy  \cref{cchi}-\cref{uvp}. Then there exists a generalised supersolution $(u,v)$ { to system \cref{KS}} in terms of {\cref{defsol}}. 
\end{theorem}
The proof of this theorem is based on a suitable regularisation and a series of priori estimates in {\cref{secreg}} leading into  a limit procedure in {\cref{SecEx}}. 
\begin{remark}[Choice of $M$]
As we see later on in the proof, a solution in terms of {\cref{defsol}} can be constructed using any function $M$ which satisfies
\begin{align}
 M{>} M_{*}{ [c]},
\end{align}
where, { for any function $\sigma\in C(\R_0^+\times \R_0^+)$,}
\begin{align}
 M_{*}{[\sigma]}=&M_{*}{[\sigma]}(n,l)\nonumber\\
 :=&{\frac{1}{4}}\frac{ab}{(a-1)(b-1)}\underset{(u,v)\in[0,n]\times[0,l]}{\max}\left\{\frac{(a-1)^2}{b^2}{\sigma}^2(u,v)\overline{u}^{a-2}\overline{v}^2+4\frac{a+b-1}{ab}\overline{u}^{a}\right\}.\label{condM}
\end{align}

\end{remark}

Our interest in the introduced concept of generalised supersolutions is supported by the following result:
\begin{theorem}[Classical solutions]\label{classol} 
{ Let 
\begin{align}
\chi\in C^1(\R_0^+\times \R_0^+).\label{chi1}\end{align}
}Let a pair $(u,v)$ be a supersolution in terms of {\cref{defsol}} { such} that for some $T\in{(0,\infty]}$
 \begin{align}
 u,v\in C([0,T)\times\overline{\Omega})\cap C^{1,2}((0,T)\times\overline{\Omega}).\label{uvsmooth} 
\end{align}
 Then $(u,v)$ is a classical solution to   \cref{KS} in ${ [0,T)}\times \overline{\Omega}$. 
\end{theorem}
The proof of {\cref{classol}} is given in {\cref{SecCl}}.

\section{Smooth regularisations for  \texorpdfstring{\cref{KS}}{}}\label{secreg}
{ Consider a sequence of functions \begin{align*}
\chi_k\in C_0^2(\R_0^+\times \R_0^+),\qquad k\in\N,                                              \end{align*}
which approximates $\chi$ on compact sets:
\begin{align*}
 \chi_k\underset{k\rightarrow\infty}{\rightarrow}\chi\qquad \text{in }C([0,n]\times[0,l])\qquad\text{for all }n,l\in\N.
\end{align*}
Set
\begin{align*}
 c_k(u,v):=\chi_k(u,v)u\qquad\text{for all }u,v\in\R_0^+.
\end{align*}
The assumptions on $\chi_k$ readily imply that
\begin{align*}
  &c_k\in C_0^2(\R_0^+\times \R_0^+), \\ 
  & c_k\underset{k\rightarrow\infty}{\rightarrow}c\qquad \text{in }C([0,n]\times[0,l])\qquad\text{for all }n,l\in\N,
\end{align*} 
and 
\begin{align}
 M_{*}[c_k](n,l)\underset{k\rightarrow\infty}{\rightarrow}M_{*}[c](n,l)\qquad\text{for all }n,l\in\N,\label{convM*}
\end{align}
with $M_*$ as defined in \cref{condM}.
Now we are ready to introduce} a family of regularisations of system \cref{KS}: 
\begin{subequations}\label{KSe}
 \begin{alignat}{3}
  &\partial_t u_k=\nabla\cdot\left(\nabla u_k-{ c_k}(u_k,v_k)\nabla v_k\right)\quad&&\text{ in }(0,T)\times\Omega,\label{ue}\\
  &\partial_t v_k=\Delta v_k-v_k+{ u_k}\quad&&\text{ in }(0,T)\times\Omega,\label{ve}\\
  &\partial_{\nu} u_k=\partial_{\nu} v_k=0\quad&&\text{ in }(0,T)\times\partial\Omega,\label{uvebc}\\
  &u_k(0,\cdot)=u_{k0},\ v_k(0,\cdot)=v_{k0}\quad&&\text{ in }\Omega.
 \end{alignat}
\end{subequations}
We choose the regularised initial data  $u_{k0}$ and $v_{k0}$ so as to satisfy
\begin{align}
&0< u_{k0},v_{k0}\in W^{1,\infty}(\Omega),\label{RegIni}\\
&\|v_{k0}\|_{L^1(\Omega)}\geq \|v_0\|_{L^1(\Omega)},\label{aproxminv0}
\end{align}
and
\begin{alignat}{3}
&u_{k0}\underset{k\rightarrow\infty}{\rightarrow}u_0&&\qquad\text{in }L^1(\Omega)\text{ and a.e. in }\Omega,\label{aproxiniu_}\\
&v_{k0}\underset{k\rightarrow\infty}{\rightarrow}v_0&&\qquad\text{in }L^1(\Omega)\text{ and a.e. in }\Omega.\label{aproxiniv_}
\end{alignat}
Classical theory  for upper-triangular systems (see, e.g., \cite{Amann1}) implies that \cref{KSe} possesses a unique global { bounded} classical solution  $(u_k,v_k)$ { provided that both solution components are a priori bounded. The latter is the case here. Indeed, we observe first that}  due to the maximum principle  both solution components are strictly positive in $\R^+_0\times \overline{\Omega}$. { Further, since $c_k$ is compactly supported, the maximum
 principle implies that $u_k$ is bounded above, which, together with the comparison principle, leads to boundedness  of  $v_k$.  Solutions to \cref{KSe} are further} studied in  \cref{basic,VarF,keyest}. 

 Throughout the rest of this section we  assume that
 \begin{align*}
  { n,l,k\in\N}
 \end{align*}
are arbitrary but fixed.  
\begin{notation}\label{Notation2}~
\begin{enumerate}[(1)]
 \item\label{Notitem} In addition to our previous conventions (see {\cref{Notation}})  we do {\bf not} indicate in an explicit way the  dependence of quantities $B_i,C_i,D,F$  (are introduced in what follows) and that of $M_*$  (as defined by \cref{condM}) upon $n$, $l$, and $M_0$ (yet another parameter which is used below).
\item   On the other hand, we stress that all estimates which we derive in this section are uniform w.r.t. $k$. In particular, all $C_i$'s are  independent of this number.

\end{enumerate}
\end{notation}
\subsection{Basic properties of    \texorpdfstring{\cref{KSe}}{}}\label{basic}
Integrating equations \cref{ue} and \cref{ve} over $\Omega$ and using the boundary conditions and partial integration we obtain the following information about the total masses: for all $t\geq0$
\begin{align}
 &\|u_k(t,\cdot)\|_{L^1(\Omega)}= \|u_{k 0}\|_{L^1(\Omega)},\label{massue}\\
 &\|v_k(t,\cdot)\|_{L^1(\Omega)}\leq\left(1-e^{-t}\right)\|u_{k 0}\|_{L^1(\Omega)}+e^{-t}\|v_{k 0}\|_{L^1(\Omega)}.\label{massve}
\end{align}
Due to \cref{massue}-\cref{massve} and a classical result based on duality (see, e.g., the proof of Lemma 5 in \cite[Appendix A]{BOTHE2010120}) we have for all \begin{align}
(r,s)\in\left[1,\frac{N+2}{N}\right)\times\left[1,\frac{N+2}{N+1}\right)\nonumber                                                                                                                                                                \end{align}
that
\begin{align}
 \left\{\left(v_k,\nabla v_k\right)\right\}_{k\in(0,1]} \quad \text{is precompact in }L^r_{loc}(\R_0^+\times\overline{\Omega})\times L^s_{loc}(\R_0^+\times\overline{\Omega}).\label{vcomp}
\end{align}
Further, using the maximum principle and the strict positivity of the Neumann heat kernel, we conclude from   \cref{aproxminv0} that $v_{k}$ can be controlled  from below in the following way: 
\begin{align}
 \inf_{(\tau,T)\times \Omega}v_{k}\geq &\inf_{(\tau,T)\times \Omega}e^{-t}e^{t\Delta}v_{k 0}\nonumber\\
 \geq &\C{\left(\tau,T\right)}\|v_{k 0}\|_{L^1(\Omega)}\nonumber\\
 \geq &\Cl{vmin}{\left(\tau,T\right)}>0\qquad\text{for all }0<\tau<T<\infty.\label{estvmin}
\end{align}
\subsection{A variational formulation for  \texorpdfstring{\cref{KSe}}{}}\label{VarF}
Let 
\begin{align*}
 a,b>2,\qquad M_0\geq0
\end{align*}
be some numbers and consider the function
\begin{align}
 &F(u,v):=\left(M_0 +\overline{u}^{a}\right)\overline{v}^{b}.\label{F}
\end{align}
We have that
\begin{align*}
 F\in C^2_b(\R_0^+\times \R_0^+),
\end{align*}
and the  derivatives of $F$ up to order two are: 
\begin{subequations}\label{FDir}
\begin{align}
&\partial_uF(u,v)=-a\overline{u}^{a-1}\overline{v}^{b},\label{Fu}\\
&\partial_vF(u,v)=-b\left(M_0 +\overline{u}^{a}\right)\overline{v}^{b-1},\label{Fv}\\
 &\partial_{uu}F(u,v)=a(a-1)\overline{u}^{a-2}\overline{v}^{b},\label{Fuu}\\
 &\partial_{uv}F(u,v)=ab\overline{u}^{a-1}\overline{v}^{b-1},\label{Fuv}\\
 &\partial_{vv}F(u,v)=b(b-1)\left(M_0+\overline{u}^{a}\right)\overline{v}^{b-2}.\label{Fvv}
\end{align}
\end{subequations}
Multiplying  \cref{ue} and \cref{ve} by $\partial_{u}F(u_k,v_k)$ and $\partial_{v}F(u_k,v_k)$, respectively, adding the results together, and using the chain rule and expressions \cref{FDir} where necessary, we compute  that
\begin{align}
 &\partial_t\left(\left(M_0 +\overline{u_k}^{a}\right)\overline{v_k}^{b}\right) \nonumber\\
 =&-\left(\nabla u_k-{ c_k}(u_k,v_k)\nabla v_k\right)\cdot\left(a(a-1)\overline{u_k}^{a-2}\overline{v_k}^{b}\nabla u_k+ab\overline{u_k}^{a-1}\overline{v_k}^{b-1}\nabla v_k\right)\nonumber\\
 &+\nabla\cdot\left(-a\overline{u_k}^{a-1}\overline{v_k}^{b}\left(\nabla u_k-{ c_k}(u_k,v_k)\nabla v_k\right)\right)\nonumber\\
 &-\nabla v_k \cdot\left(ab\overline{u_k}^{a-1}\overline{v_k}^{b-1}\nabla u_k+b(b-1)\left(M_0+\overline{u_k}^{a}\right)\overline{v_k}^{b-2}\nabla v_k\right)\nonumber\\
 &+\nabla\cdot\left(-b\left(M_0 +\overline{u_k}^{a}\right)\overline{v_k}^{b-1}\nabla v_k\right)\nonumber\\
 &-b\left(M_0 +\overline{u_k}^{a}\right)\overline{v_k}^{b-1}\left(-v_k+{ u_k}\right)
 \nonumber\\
 =&-\left(a(a-1)\overline{u_k}^{a-2}\overline{v_k}^{b}\left|\nabla u_k\right|^2\right.\nonumber\\
 &\phantom{aaa}\left.+\left(2ab\overline{u_k}^{a-1}\overline{v_k}^{b-1}-a(a-1){ c_k}(u_k,v_k)\overline{u_k}^{a-2}\overline{v_k}^{b}\right)\nabla u_k\cdot\nabla v_k\right.\nonumber\\
 &\phantom{aaa}\left.+\left(b(b-1)\left(M_0+\overline{u_k}^{a}\right)\overline{v_k}^{b-2}-ab{ c_k}(u_k,v_k)\overline{u_k}^{a-1}\overline{v_k}^{b-1}\right)\left|\nabla v_k\right|^2\right)\nonumber\\
 &+\nabla\cdot\left(-a\overline{u_k}^{a-1}\overline{v_k}^{b}\nabla u_k+\left(-b\left(M_0 +\overline{u_k}^{a}\right)\overline{v_k}^{b-1}+a{ c_k}(u_k,v_k)\overline{u_k}^{a-1}\overline{v_k}^{b}\right)\nabla v_k\right)\nonumber\\
 &-b\left(M_0 +\overline{u_k}^{a}\right)\overline{v_k}^{b-1}\left(-v_k+{ u_k}\right)\nonumber\\
 =&-4\left(\frac{a-1}{a}\left|\nabla \left(\overline{u_k}^{\frac{a}{2}}\overline{v_k}^{\frac{b}{2}}\right)-\overline{u_k}^{\frac{a}{2}}\nabla \overline{v_k}^{\frac{b}{2}}\right|^2\right.\nonumber\\
 &\phantom{aaaa}\left.+\frac{1}{b}\overline{u_k}^{\frac{a}{2}-1}\left(2b\overline{u_k}-(a-1){ c_k}(u_k,v_k)\overline{v_k}\right)\left(\nabla \left(\overline{u_k}^{\frac{a}{2}}\overline{v_k}^{\frac{b}{2}}\right)-\overline{u_k}^{\frac{a}{2}}\nabla \overline{v_k}^{\frac{b}{2}}\right)\cdot\nabla \overline{v_k}^{\frac{b}{2}}\right.\nonumber\\
 &\phantom{aaaa}\left.+\frac{1}{b}\left((b-1)\left(M_0+\overline{u_k}^{a}\right)-a{ c_k}(u_k,v_k)\overline{u_k}^{a-1}\overline{v_k}\right)\left|\nabla \overline{v_k}^{\frac{b}{2}}\right|^2\right)\nonumber\\
 &+\nabla\cdot\left(\left(\nabla\left( \overline{u_k}^{a}\overline{v_k}^{b}\right)-\overline{u_k}^{a}\nabla \overline{v_k}^{b}\right)+\frac{1}{b}\left(b\left(M_0 +\overline{u_k}^{a}\right)-a{ c_k}(u_k,v_k)\overline{u_k}^{a-1}\overline{v_k}\right)\nabla \overline{v_k}^{b}\right)\nonumber\\
 &-b\left(M_0 +\overline{u_k}^{a}\right)\overline{v_k}^{b-1}\left(-v_k+{ u_k}\right).
\label{KEYab_}
\end{align}
The first big bracket on the  right-hand side of  \cref{KEYab_} is a quadratic form with respect to 
\begin{align*}
 \nabla \left(\overline{u_k}^{\frac{a}{2}}\overline{v_k}^{\frac{b}{2}}\right)-\overline{u_k}^{\frac{a}{2}}\nabla \overline{v_k}^{\frac{b}{2}}\qquad\text{and}\qquad \nabla \overline{v_k}^{\frac{b}{2}}.
\end{align*}
Its discriminant $D{ [c_k]}=D{ [c_k]}(u,v)$ satisfies
\begin{align}
 &D{[\sigma]}(u,v)\nonumber\\
 =&\frac{1}{b^2}\overline{u}^{a-2}\left(2b\overline{u}-(a-1){\sigma}(u,v)\overline{v}\right)^2-4\frac{a-1}{ab}\left((b-1)\left(M_0+\overline{u}^{a}\right)-a{\sigma}(u,v)\overline{u}^{a-1}\overline{v}\right)\nonumber\\
 =&\frac{(a-1)^2}{b^2}{\sigma}^2(u,v)\overline{u}^{a-2}\overline{v}^2+4\frac{a+b-1}{ab}\overline{u}^{a}-{ 4}\frac{(a-1)(b-1)}{ab}M_0\nonumber\\
 \leq&{ 4}\frac{(a-1)(b-1)}{ab}(M_{*}{ [\sigma]}-M_0),\nonumber
\end{align}
with $M_{*}$ is as defined in \cref{condM}. Thus, \cref{KEYab_} takes the form
\begin{align}
 &\partial_tF(u_k,v_k)\nonumber\\
 =&-4\frac{a-1}{a}\left|\nabla \left(\overline{u_k}^{\frac{a}{2}}\overline{v_k}^{\frac{b}{2}}\right)+\Cl[B]{B1}{[c_k]}(u_k,v_k)\nabla \overline{v_k}^{\frac{b}{2}}\right|^2\nonumber\\
 &-\left(4\frac{b-1}{b}(M_0-M_{*}{ [c_k]})+\Cl[B]{B2}{[c_k]}(u_k,v_k)\right)\left|\nabla \overline{v_k}^{\frac{b}{2}}\right|^2\nonumber\\
 &+\nabla\cdot\left(\nabla\left(\overline{u_k}^{a}\overline{v_k}^{b} \right)+\Cl[B]{B4}{[c_k]}(u_k,v_k)\nabla \overline{v_k}^{b}\right)+{\Cl[B]{Bnew}}(u_k,v_k)-{\Cl[B]{B5}(u_k,v_k)}{ u_k},
\label{KEYab}
\end{align}
where, in order to simplify the exposition, we have introduced
\begin{align*}
 &\Cr{B1}{ [\sigma]}(u,v):=-\overline{u}^{\frac{a}{2}}+\frac{a}{2(a-1)}\frac{1}{b}\overline{u}^{\frac{a}{2}-1}\left(2b\overline{u}-(a-1){ \sigma}(u,v)\overline{v}\right),\\
 &\Cr{B2}{ [\sigma]}(u,v):=-\frac{a}{a-1}D{ [\sigma]}(u,v)-4\frac{b-1}{b}(M_0-M_{*}{ [\sigma]}),\\
 &\Cr{B4}{ [\sigma]}(u,v):=-\overline{u}^{a}+\frac{1}{b}\left(b\left(M_0 +\overline{u}^{a}\right)-a{ \sigma}(u,v)\overline{u}^{a-1}\overline{v}\right),\\
 &{\Cr{Bnew}(u,v):=b\left(M_0 +\overline{u}^{a}\right)\overline{v}^{b-1}v,}\\
 &\Cr{B5}(u,v):=b\left(M_0 +\overline{u}^{a}\right)\overline{v}^{b-1}.
\end{align*}
Observe that
\begin{alignat}{5}
 &{ B_i[c_k]\in C_b(\R_0^+\times\R_0^+),\qquad }&&{\|B_i[c_k]\|_{C_b(\R_0^+\times\R_0^+)}\leq{\Cl{CB}}\qquad}&&{\text{for }i\in\{ 1,2,3\},}\label{Bbound1}\\
 &B_i\in C_b(\R_0^+\times\R_0^+),\qquad &&\|B_i\|_{C_b(\R_0^+\times\R_0^+)}\leq{\Cr{CB}}\qquad&&\text{for }i\in\{{ 4,5}\},\label{Bbound2}\\
 &\Cr{B2}{[c_k]},\Cr{B5}\geq0,&&&&\label{Bpos}
\end{alignat}
and
\begin{align}
 &{ B_i[c_k]\underset{k\rightarrow\infty}{\rightarrow}B_i[c]\qquad\text{in }C_b(\R_0^+\times\R_0^+)\qquad\text{for }i\in\{ 1,2,3\}.}\label{Bconv}
\end{align}
Multiplying \cref{KEYab} by an arbitrary function $\psi\in C^1_0(\R_0^+)$ and integrating by parts w.r.t. $t$ yields for all $x\in\Omega$ that
\begin{align}
 &-\int_0^{\infty}F(u_k,v_k)\partial_t\psi\,dt-F(u_{k0},v_{k0})\psi(0)\nonumber\\
 =&-\int_0^{\infty}4\frac{a-1}{a}\left|\nabla \left(\overline{u_k}^{\frac{a}{2}}\overline{v_k}^{\frac{b}{2}}\right)+\Cr{B1}{[c_k]}(u_k,v_k)\nabla \overline{v_k}^{\frac{b}{2}}\right|^2\psi\,dt\nonumber\\
 &-\int_0^{\infty}\left(4\frac{b-1}{b}(M_0-M_{*}{ [c_k]})+\Cr{B2}{[c_k]}(u_k,v_k)\right)\left|\nabla \overline{v_k}^{\frac{b}{2}}\right|^2\psi\,dt\nonumber\\
 &+\nabla\cdot\int_0^{\infty}\psi\left(\nabla\left(\overline{u_k}^{a}\overline{v_k}^{b} \right)+\Cr{B4}{[c_k]}(u_k,v_k)\nabla \overline{v_k}^{b}\right)\,dt\nonumber\\
 &+\int_0^{\infty}{\Cr{Bnew}}(u_k,v_k)\psi-{\Cr{B5}}(u_k,v_k){ u_k}\psi\,dt.
\label{KEYabt}
\end{align}
On the other hand, multiplying \cref{KEYab} by an arbitrary function $\varphi\in C^1(\Omega)$ and integrating by parts w.r.t. $x$ using the boundary conditions where necessary yields that
\begin{align}
 \int_{\Omega}\partial_tF(u_k,v_k)\varphi\,dx
 =&-\int_{\Omega}4\frac{a-1}{a}\left|\nabla \left(\overline{u_k}^{\frac{a}{2}}\overline{v_k}^{\frac{b}{2}}\right)+\Cr{B1}{[c_k]}(u_k,v_k)\nabla \overline{v_k}^{\frac{b}{2}}\right|^2{\varphi}\,dx\nonumber\\
 &-\int_{\Omega}\left(4\frac{b-1}{b}(M_0-M_{*}{ [c_k]})+\Cr{B2}{[c_k]}(u_k,v_k)\right)\left|\nabla \overline{v_k}^{\frac{b}{2}}\right|^2{\varphi}\,dx\nonumber\\
 &-\int_{\Omega}\left(\nabla\left(\overline{u_k}^{a}\overline{v_k}^{b} \right)+\Cr{B4}{[c_k]}(u_k,v_k)\nabla \overline{v_k}^{b}\right)\cdot\nabla\varphi\,dx\nonumber\\
 &+\int_{\Omega}{\Cr{Bnew}}(u_k,v_k)\varphi-{\Cr{B5}}(u_k,v_k){ u_k}\varphi\,dx.
\label{KEYabx}
\end{align}
Finally, if we multiply \cref{KEYab} by the product $\psi\varphi$ and integrate by parts w.r.t. $t$ and $x$, then we arrive at the following variational reformulation:
\begin{align}
 &-\int_0^{\infty}\partial_t\psi\int_{\Omega}F(u_k,v_k)\varphi\,dxdt-\psi(0)\int_{\Omega}F(u_{k0},v_{k0}){\varphi}\,dx\nonumber\\
 =&-\int_0^{\infty}\psi\int_{\Omega}4\frac{a-1}{a}\left|\nabla \left(\overline{u_k}^{\frac{a}{2}}\overline{v_k}^{\frac{b}{2}}\right)+\Cr{B1}{[c_k]}(u_k,v_k)\nabla \overline{v_k}^{\frac{b}{2}}\right|^2{\varphi}\,dxds\nonumber\\
 &-\int_0^{\infty}\psi\int_{\Omega}\left(4\frac{b-1}{b}(M_0-M_{*}{ [c_k]})+\Cr{B2}{[c_k]}(u_k,v_k)\right)\left|\nabla \overline{v_k}^{\frac{b}{2}}\right|^2{\varphi}\,dxds\nonumber\\
 &-\int_0^{\infty}\psi\int_{\Omega}\left(\nabla\left(\overline{u_k}^{a}\overline{v_k}^{b} \right)+\Cr{B4}{[c_k]}(u_k,v_k)\nabla \overline{v_k}^{b}\right)\cdot\nabla\varphi\,dxds\nonumber\\
 &+\int_0^{\infty}\psi\int_{\Omega}{\Cr{Bnew}}(u_k,v_k)\varphi-{\Cr{B5}}(u_k,v_k){ u_k}\varphi\,dxds.
\label{superue}
\end{align}
For equation \cref{ve} a standard procedure yields the following reformulations:
for all ${\varphi\in C^1(\Omega)}$ and ${\psi\in C^1_0(\R_0^+)}$ it holds that
\begin{align}
 -\int_0^{\infty}v_k\partial_t\psi \,ds-\psi(0)v_{k0}=\nabla\cdot\int_0^{\infty}\psi\nabla v_k\,dt+\int_0^{\infty}\left(-v_k+{ u_k}\right)\psi\,dt.\label{supervte}
\end{align}
and
\begin{align}
 -\int_0^{\infty}\partial_t\psi\int_{\Omega}v_k\varphi\,dxds-\psi(0)\int_{\Omega}v_{k0}\varphi\,dx
 =&\int_0^{\infty}\psi\int_{\Omega}-\nabla v_k\cdot\nabla\varphi+\left(-v_k+{ u_k}\right)\varphi\,dxds.\label{superve}
\end{align}

\subsection{Further uniform \texorpdfstring{{ (w.r.t. $k$)}}{} estimates for \texorpdfstring{\cref{KSe}}{}}\label{keyest}
Choosing $\varphi\equiv 1$ and 
\begin{align*}
M_0={ M_{*}[c_k]+1}
\end{align*}
in \cref{KEYabx}  and using \cref{Bbound1}-\cref{Bpos} yields 
that 
\begin{align}
 &\frac{d}{dt}\int_{\Omega}{ \left(M_{*}[c_k]+1 +\overline{u}_k^{a}\right)\overline{v}_k^{b}}\,dx\nonumber\\
 =&-\int_{\Omega}4\frac{a-1}{a}\left|\nabla \left(\overline{u_k}^{\frac{a}{2}}\overline{v_k}^{\frac{b}{2}}\right)+\Cr{B1}{[c_k]}(u_k,v_k)\nabla \overline{v_k}^{\frac{b}{2}}\right|^2\,dx\nonumber\\
 &-\int_{\Omega}\left({ 4\frac{b-1}{b}}+\Cr{B2}{[c_k]}(u_k,v_k)\right)\left|\nabla \overline{v_k}^{\frac{b}{2}}\right|^2\,dx\nonumber\\
 &+\int_{\Omega}{\Cr{Bnew}}(u_k,v_k)-{\Cr{B5}}(u_k,v_k){ u_k}\,dx\nonumber\\
 \leq&-\int_{\Omega}4\frac{a-1}{a}\left|\nabla \left(\overline{u_k}^{\frac{a}{2}}\overline{v_k}^{\frac{b}{2}}\right)+\Cr{B1}{[c_k]}(u_k,v_k)\nabla \overline{v_k}^{\frac{b}{2}}\right|^2{+4\frac{b-1}{b}}\left|\nabla \overline{v_k}^{\frac{b}{2}}\right|^2\,dx+\Cr{CB}|\Omega|.
\label{KeyApr}
\end{align}
 Integrating the differential inequality \cref{KeyApr}  over $(0,T)$ for any $T>0$ we conclude that
\begin{align}
 &\int_{\Omega}{ \left(M_{*}[c_k]+1 +\overline{u}_k^{a}\right)\overline{v}_k^{b}}(T)\,dx\nonumber\\
 &+ \int_0^T\int_{\Omega}4\frac{a-1}{a}\left|\nabla \left(\overline{u_k}^{\frac{a}{2}}\overline{v_k}^{\frac{b}{2}}\right)+\Cr{B1}{[c_k]}(u_k,v_k)\nabla \overline{v_k}^{\frac{b}{2}}\right|^2{+4\frac{b-1}{b}}\left|\nabla \overline{v_k}^{\frac{b}{2}}\right|^2\,dxds\nonumber\\
 \leq& \int_{\Omega}{ \left(M_{*}[c_k]+1 +\overline{u}_{k0}^{a}\right)\overline{v}_{k0}^{b}}\,dx+\C(T).\label{est1}
\end{align}
{ Since $M_{*}[c_k]$ is nonnegative and uniformly bounded w.r.t. $k$ (recall \cref{convM*}),} 
the integral inequality \cref{est1} directly implies 
\begin{align}
 &\int_0^T\int_{\Omega}\left|\nabla \left(\overline{u_k}^{\frac{a}{2}}\overline{v_k}^{\frac{b}{2}}\right)+\Cr{B1}{[c_k]}(u_k,v_k)\nabla \overline{v_k}^{\frac{b}{2}}\right|^2\,dxds\leq \Cl{C1new}(T),\label{estnew}\\
 &{\int_0^T}\int_{\Omega}\left|\nabla \overline{v_k}^{\frac{b}{2}}\right|^2\,dxds\leq\Cr{C1new}(T).\label{est12new}
\end{align}
Combining \cref{estnew}-\cref{est12new} and using the fact that { $\Cr{B1}{[c_k]}$ is uniformly bounded w.r.t. $k$}, we also get
\begin{align}
 &{\int_0^T}\int_{\Omega}\left|\nabla \left(\overline{u_k}^{\frac{a}{2}}\overline{v_k}^{\frac{b}{2}}\right)\right|^2\,dxds\leq\Cl{C2new}(T),\label{est11new}
\end{align}
which due to the boundedness of the cut-offs immediately implies that
\begin{align}
 &{\int_0^T}\left\| \overline{u_k}^{\frac{a}{2}}\overline{v_k}^{\frac{b}{2}}\right\|_{H^1(\Omega)}^2\,ds\leq\Cl{C2new1}(T),\label{est11new1}
\end{align}
Taking $M_0=0$ in \cref{KEYabx} we obtain { thanks to  \cref{estnew}-\cref{est11new}} and the boundedness of $F$ and $B_i$'s that for a $p>N$
\begin{align}
 &\int_0^T\left\|\partial_t\left(\overline{u_k}^{a}\overline{v_k}^{b}\right)\right\|_{(W^{1,p}(\Omega))^*}\leq \Cl{C3}({ T}).\label{estprt}
\end{align}
Further, using \cref{estnew}-\cref{est11new}  and the boundedness of $F$ and $B_i$'s, we conclude from  \cref{KEYabt}
that
\begin{align}
 \left\|\nabla\cdot\int_0^{\infty}\psi\left(\nabla\left(\overline{u_k}^{a}\overline{v_k}^{b} \right)+\Cr{B4}{[c_k]}(u_k,v_k)\nabla \overline{v_k}^{b}\right)\,dt\right\|_{L^1(\Omega)}\leq {\C(\psi)}.\label{estintnabla_}
\end{align}
Similarly, using \cref{aproxiniv_} and \cref{massue}-\cref{massve} we deduce from \cref{supervte} that
\begin{align}
 \left\|\nabla\cdot\int_0^{\infty}\psi\nabla v_{k}\,dt\right\|_{L^1(\Omega)}\leq \C(\psi).\label{estintDelta_}
\end{align}
\section{Construction of a generalised supersolution to \texorpdfstring{\cref{KS}}{}: proof of {\texorpdfstring{\cref{TheoEx}}{}}}\label{SecEx}
\begin{notation}
 Below $k_m$ always denotes a suitably chosen subsequence, i.e.,   such that the stated convergence holds. It may thus become 'thinner' from one instance to another.
%
\end{notation}
\subsection{Convergence in the \texorpdfstring{$v$-}{}component}\label{SubSecv}

To begin with, we use  properties established in {\cref{basic}} and  conclude that there exists a  measurable function $v:[0,\infty)\times \overline{\Omega}\rightarrow\R^+$ such that:\\
due to \cref{vcomp},
\begin{alignat}{3}
 &v_{k_m}\underset{m\rightarrow\infty}{\rightarrow}v&&\qquad\text{in }L^r_{loc}(\R_0^+\times\overline{\Omega})\text{ and a.e. in }\R^+\times\Omega,\label{convv}\\
 &\nabla v_{k_m}\underset{m\rightarrow\infty}{\rightarrow}\nabla v&&\qquad\text{in }L^s_{loc}(\R_0^+\times\overline{\Omega})\text{ and a.e. in }\R^+\times\Omega;\label{convnv}
\end{alignat}
due to \cref{estvmin} and \cref{convv}, for all $0<\tau<T<\infty$
\begin{align}
 \underset{(\tau,T)\times \Omega}{\ess\inf}v\geq \Cr{vmin}(\tau,T)>0,\label{vmin}
\end{align}
so that $v>0$ a.e. in $\R^+\times\Omega$. 
\subsection{Almost everywhere convergence in the \texorpdfstring{$u$-}{}component}\label{SubSecue}
Next, we address the a.e. convergence of the $u$-component. We deduce that:\\
{ due to \cref{est11new1}} for $a=b=6$ and \cref{estprt} for $a=b=3$, the Lions-Aubin lemma \cite[Corollary 4]{Simon},  and the continuity of the root function, 
\begin{align}
&{\text{for each }n,l\in\N\text{ every  subsequence of }}\left(n-u_{k_m}\right)_+\left(l-v_{k_m}\right)_+\text{ has a subsequence}\nonumber\\
&\text{that converges  a.e. in }\R^+\times\Omega;\label{convpr1_}
\end{align}
due to \cref{convv} and \cref{convpr1_}, 
\begin{align}
 &{\text{for each }n,l\in\N\text{  every  subsequence of }}\left(n-u_{k_m}\right)_+\text{ has a subsequence}\nonumber\\
&\text{that converges a.e. in }\{v<l\};\label{connmu1}
\end{align}
due to \cref{connmu1}, $v\in L^1_{loc}(\R_0^+\times\overline{\Omega})$, and   a diagonal argument,
\begin{align}
&\text{for each }n\in\N\text{  every  subsequence of }\left(n-u_{k_m}\right)_+\text{ has a subsequence}\nonumber\\
&\text{that converges a.e. in }\
 \R^+\times\Omega;\label{connmu2}
\end{align}
due to \cref{connmu2} and a diagonal argument,
\begin{align}
 \text{for each }n\in\N\text{ sequence } \left(n-u_{k_m}\right)_+\text{ converges a.e. in }\R^+\times\Omega;\label{connmu3_}
\end{align}
due to \cref{connmu3_},
\begin{align}
 \text{for each }n\in\N\text{ sequence } \min\{u_{k_m},n\}\text{ converges to }{\eta}_n\in[0,n]
 \text{ a.e. in }\R^+\times\Omega\label{connmu3}
\end{align}
for some measurable ${\eta}_n$;\\
due to \cref{connmu3} and $n\mapsto \min\{u,n\}$ increasing for all $u\in\R$,
\begin{align}
 {\eta}_n\underset{n\rightarrow\infty}{\nearrow}{ u}\in[0,\infty]\text{ a.e. in }\R^+\times\Omega\label{etaincr}
\end{align}
for some measurable $u$;\\
due to \cref{connmu3},
\begin{align}
 \text{for each }n\in\N\text{ sequence } u_{k_m}\text{ converges to }{\eta}_n\text{ a.e. in }\{{\eta}_n<n\};\label{connmu4}
\end{align}
due to \cref{connmu4,etaincr},
\begin{align}
 {\text{for each }n\in\N\text{ sequence }}u_{k_m}\text{ converges to }{\eta}_n{ =u}\ \ \text{a.e. in }{\{{ u}<n\}\subset\{{\eta}_n<n\}};\label{connmu5}
\end{align}
due to \cref{aproxiniu_}, \cref{massue}, \cref{connmu3}, \cref{etaincr},  $\min\{u,n\}\in[0,u]$, and the  monotone convergence theorem,  { $u\in L^{\infty}(\R^+;L^1(\Omega))$ and}
\begin{align}
 \|{ u}\|_{L^{\infty}(\R^+;L^1(\Omega))}\leq \|u_0\|_{L^1(\Omega)};\label{xiL1}
\end{align}
 due to $u\in L^{\infty}(\R^+;L^1(\Omega))$, \cref{connmu5}, and a diagonal argument,
\begin{align}
 &u_{k_m}
  \underset{m\rightarrow\infty}{\rightarrow}u\qquad\text{ a.e. in }\R^+\times\Omega.\label{convuae}
\end{align}
\subsection{Convergence in the variational reformulations}
In this sequel we prove that $u$ and $v$ satisfy variational reformulations \cref{superu} and \cref{prweakbc}-\cref{vweakbc}. 
 From now on we fix arbitrary constants
\begin{align*}
 a,b>2
\end{align*}
and a function 
\begin{align}
 &M:\N\times\N\rightarrow \R^+_0,\nonumber\\
 &{ M{ >} M_*{[c]}},\label{Ml}
\end{align}
with $M_*$ as defined in \cref{condM}. Thus, these parameters satisfy the requirements of {\cref{defsol}}. Further, throughout this section we assume $$n,l\in\N$$ to be arbitrary but fixed and make use of {\cref{Notation}} as well as {\cref{Notation2}\eqref{Notitem}}.

{ Consider an arbitrary sequence of functions
 \begin{align*}
&G_{{ m}}\in C_b(\R^+_0\times \R^+_0),
\end{align*}
such that for some $G\in C_b(\R^+_0\times \R^+_0)$
\begin{align*}
 G_{{ m}}\underset{m\rightarrow\infty}{\rightarrow}G\qquad\text{in } C_b(\R^+_0\times \R^+_0).
\end{align*}
}Using the convergence results from \cref{SubSecv,SubSecue}, we deduce that:\\
due to \cref{convv} and \cref{convuae},
\begin{alignat}{3}
 &G_{{ m}}\left(u_{k_m},v_{k_m}\right)\underset{m\rightarrow\infty}{\rightarrow}G(u,v)&&\qquad\text{ a.e. in }\R^+\times\Omega;\label{convG}
\end{alignat}
due to \cref{convG} and the dominated convergence theorem,
\begin{align}
 &G_{{ m}}\left(u_{k_m},v_{k_m}\right)\underset{m\rightarrow\infty}{\rightarrow}G(u,v)\qquad\text{ in }L^p_{loc}(\R^+_0\times\Omega)\text{ for all }p>1;\label{convGLp}
\end{align}
due to \cref{convG},
\begin{align}
 G_{{ m}}(u_{k_m},v_{k_m})\underset{m\rightarrow\infty}{\rightarrow}{ G}(u,v)\qquad\text{ in }L^1_{loc}(\R^+_0\times\Omega);\label{convFL1}
\end{align}
due to \cref{convv}-\cref{convnv} and the chain rule,
\begin{align}
 \nabla \overline{v_{k_m}}^{\frac{b}{2}}\underset{m\rightarrow\infty}{\rightarrow}\nabla \overline{v}^{\frac{b}{2}}\qquad \text{in }L^s_{loc}(\R^+_0\times\Omega)\text{ and a.e. in }\R^+\times\Omega;\label{convnbv}
\end{align}
{ due to \cref{est12new}} and  \cref{convnbv}, 
\begin{align}
 \nabla \overline{v_{k_m}}^{\frac{b}{2}}\underset{m\rightarrow\infty}{\rightharpoonup}\nabla \overline{v}^{\frac{b}{2}}\qquad \text{in }L^2_{loc}(\R^+_0\times\Omega);\label{convweak}
\end{align}
due to \cref{convG} and \cref{convnbv},
\begin{align}
 G_{{ m}}\left(u_{k_m},v_{k_m}\right)\nabla \overline{v_{k_m}}^{\frac{b}{2}}\underset{m\rightarrow\infty}{\rightarrow}G\left(u,v\right)\nabla \overline{v}^{\frac{b}{2}}\qquad\text{ a.e. in }\R^+\times\Omega;\label{convGvae}
\end{align}
 due to \cref{est12new}, \cref{convGvae}, and the Lions lemma \cite[Lemma 1.3]{Lions} 
\begin{align}
  G_{{ m}}\left(u_{k_m},v_{k_m}\right)\nabla \overline{v_{k_m}}^{\frac{b}{2}}\underset{m\rightarrow\infty}{\rightharpoonup}G\left(u,v\right)\nabla \overline{v}^{\frac{b}{2}}\qquad \text{in }L^2_{loc}(\R^+_0\times\Omega);\label{convGvweak}
\end{align}
 due to \cref{est11new}, \cref{convG}, and the Banach-Alaoglu theorem,
\begin{align}
 \nabla \left(\overline{u_{k_m}}^{\frac{a}{2}}\overline{v_{k_m}}^{\frac{b}{2}}\right)\underset{m\rightarrow\infty}{\rightharpoonup}\nabla \left(\overline{{ u}}^{\frac{a}{2}}\overline{{ v}}^{\frac{b}{2}}\right)\qquad \text{in }L^2_{loc}(\R^+_0\times\Omega).\label{convprw}
\end{align}
Let 
\begin{align}
 0\leq\varphi\in C^1(\overline{\Omega}),\qquad 0\leq\psi\in C^1_0(\R_0^+)\nonumber
\end{align}
be arbitrary. Then:\\
{ due to \eqref{Bconv} and} \cref{convGvweak}-\cref{convprw}, 
\begin{align}
&\nabla \left(\overline{u_{k_m}}^{\frac{a}{2}}\overline{v_{k_m}}^{\frac{b}{2}}\right)+\Cr{B1}{[c_{k_m}]}(u_{k_m},v_{k_m})\nabla \overline{v_{k_m}}^{\frac{b}{2}}\underset{m\rightarrow\infty}{\rightharpoonup}\nabla \left(\overline{u}^{\frac{a}{2}}\overline{v}^{\frac{b}{2}}\right)+\Cr{B1}{[c]}(u,v)\nabla \overline{v}^{\frac{b}{2}},\label{convweak1}\\
 &\nabla\left(\overline{u_{k_m}}^{a}\overline{v_{k_m}}^{b} \right)+\Cr{B4}{[c_{k_m}]}(u_{k_m},v_{k_m})\nabla \overline{v_{k_m}}^{b}\underset{m\rightarrow\infty}{\rightharpoonup}\nabla\left(\overline{u}^{a}\overline{v}^{b} \right)+\Cr{B4}{[c]}(u,v)\nabla \overline{v}^{b}\label{convweak2}\\
 &\text{in }L^2_{loc}(\R^+_0\times\Omega);\nonumber
\end{align}
due to \cref{convweak1},
\begin{align}
 \underset{m\rightarrow\infty}{\lim\inf}\ &\int_0^{\infty}\psi\int_{\Omega}4\frac{a-1}{a}\left|\nabla \left(\overline{u_{k_m}}^{\frac{a}{2}}\overline{v_{k_m}}^{\frac{b}{2}}\right)+\Cr{B1}{[c_{k_m}]}(u_{k_m},v_{k_m})\nabla \overline{v_{k_m}}^{\frac{b}{2}}\right|^2\varphi\,dxds\nonumber\\
 \geq &\int_0^{\infty}\psi\int_{\Omega}4\frac{a-1}{a}\left|\nabla \left(\overline{u}^{\frac{a}{2}}\overline{v}^{\frac{b}{2}}\right)+\Cr{B1}{[c]}(u,v)\nabla \overline{v}^{\frac{b}{2}}\right|^2\varphi\,dxds;\label{conv9w}
\end{align}
 { due to \cref{convM*}}, \cref{Bpos}{-\eqref{Bconv}}, { \cref{Ml}}, \cref{convGvae}, and Fatou's lemma,  taking $M_0:=M$
\begin{align}
 \underset{m\rightarrow\infty}{\lim\inf}\ &\int_0^{\infty}\psi\int_{\Omega}\left(4\frac{b-1}{b}(M-M_{*}{ [c_{k_m}]})+\Cr{B2}{[c_{k_m}]}(u_{k_m},v_{k_m})\right)\left|\nabla \overline{v_{k_m}}^{\frac{b}{2}}\right|^2\varphi\,dxds\nonumber\\
 \geq &\int_0^{\infty}\psi\int_{\Omega}\left(4\frac{b-1}{b}(M-M_{*}{ [c]})+\Cr{B2}{[c]}(u,v)\right)\left|\nabla \overline{v}^{\frac{b}{2}}\right|^2\varphi\,dxds\label{conv11w}
\end{align}
 { due to \cref{Bbound2}-\cref{Bpos}}, \cref{convv}, \cref{convuae}, \cref{convG}, and Fatou's lemma,
\begin{align}
 \underset{m\rightarrow\infty}{\lim\sup}\ &\int_0^{\infty}\psi\int_{\Omega}{\Cr{Bnew}}(u_{k_m},v_{k_m})\varphi-{\Cr{B5}}(u_{k_m},v_{k_m}){ u_{k_m}}\varphi\,dxds\nonumber\\
 \leq& \int_0^{\infty}\psi\int_{\Omega}{\Cr{Bnew}}(u,v)\varphi-{\Cr{B5}}(u,v){ u}\varphi\,dxds;\label{conv10w}
\end{align}
due to \cref{aproxiniu_}-\cref{aproxiniv_} and the dominated convergence theorem,
\begin{align}
 F\left(u_{k_m0},v_{k_m0}\right)\underset{m\rightarrow\infty}{\rightarrow}F(u_0,v_0)\qquad\text{ in }L^1_{loc}(\R^+_0\times\Omega).\label{convini}
\end{align}
Combining \cref{convFL1} and  \cref{conv9w}-\cref{convini}, we can pass to the limit superior on both sides of \cref{superue} and thus finally obtain that  $(u,v)$ satisfies \cref{superu}.

Let now  
\begin{align}
 \varphi\in C^1(\overline{\Omega}),\qquad \psi\in C^1_0(\R_0^+),\nonumber
\end{align}
i.e., not necessarily nonnegative. Then:\\
due to {\eqref{Bconv},}  \cref{estintnabla_}, \cref{convweak2}, and the Banach-Alaoglu theorem, 
\begin{align}
 &\nabla\cdot\int_0^{\infty}\psi\left(\nabla\left(\overline{u_{k_m}}^{a}\overline{v_{k_m}}^{b} \right)+\Cr{B4}{[c_{k_m}]}(u_{k_m},v_{k_m})\nabla \overline{v_{k_m}}^{b}\right)\,dt\nonumber\\
 \underset{m\rightarrow\infty}{\overset{*}{\rightharpoonup}}&\nabla\cdot\int_0^{\infty}\psi\left(\nabla\left(\overline{u}^{a}\overline{v}^{b} \right)+\Cr{B4}{[c]}(u,v)\nabla \overline{v}^{b}\right)\,dt\qquad\text{in }{\cal M}(\overline{\Omega}),\label{convbcpr}
\end{align}
and 
\begin{align*}
\int_0^{\infty}\psi\left(\nabla\left(\overline{u}^{a}\overline{v}^{b} \right)+\Cr{B4}{[c]}(u,v)\nabla \overline{v}^{b}\right)\,dt\in {\cal DM}^2(\Omega);   
\end{align*}
due to \cref{contitr}, \cref{uvebc}, \cref{convweak2}, and \cref{convbcpr}, $(u,v)$ satisfies  \cref{prweakbc};\\
due to \cref{estintDelta_}, \cref{convnv}, and the Banach-Alaoglu theorem, 
\begin{align}
 \nabla\cdot\int_0^{\infty}\psi \nabla v_{k_m}\,dt\underset{m\rightarrow\infty}{\overset{*}{\rightharpoonup}}\nabla\cdot\int_0^{\infty}\psi\nabla v\,dt\qquad\text{in }{\cal M}(\overline{\Omega}),\label{convDelta}
\end{align}
and 
\begin{align*}
\int_0^{\infty}\psi \nabla v\,dt\in {\cal DM}^{s}(\Omega);                                                              \end{align*}
due to \cref{contitr}, \cref{uvebc}, \cref{convnv}, and \cref{convDelta},  $v$ satisfies  \cref{vweakbc}.
\subsection{Convergence to \texorpdfstring{$\mu$}{}}
Finally, we establish \cref{superv}-\cref{massu0} for some element $\mu\in L^{\infty}_{w-*}(\R^+;{\cal M}_+(\overline{\Omega}))$. It holds that:\\
due to \cref{aproxiniu_} and  \cref{massue},
\begin{align}
 \int_{\Omega}u_{k_m}\,dx
 \phantom{\ }\equiv\phantom{\    }&\int_{\Omega}u_{k_m0}\,dx\nonumber\\
  \underset{m\rightarrow\infty}{\rightarrow}&\int_{\Omega}u_0\,dx.\label{conv9}
\end{align}
Now let
\begin{align*}
 { 0\leq}\varphi\in C(\overline{\Omega}).
\end{align*}
Then:\\
{ due to \cref{aproxiniu_}, \cref{massue},} 
and the Banach-Alaoglu theorem (compare also {\cref{SecFSp}}), there exists some $\mu\in L^{\infty}_{w-*}(\R^+;{\cal M}(\overline{\Omega}))$, such that 
\begin{align}
&\int_{\Omega} u_{k_m}\varphi\,dx
 \underset{m\rightarrow\infty}{\overset{*}{\rightharpoonup}}\int_{\Omega} u\varphi\,dx+\int_{\overline{\Omega}}\varphi\ (d\mu(\cdot))(x)\qquad\text{in }L^{\infty}(\R^+);\label{conv06}
\end{align}
due to \cref{convuae} and Fatou's lemma, 
\begin{align}
 \underset{m\rightarrow\infty}{\lim}\int_{\Omega} u_{k_m}\varphi\,dx
 \geq\int_{\Omega}u\varphi\,dx\qquad\text{a.e. in }\R^+;\label{conv7}
\end{align}
due to \cref{conv06}  and  \cref{conv7}, 
\begin{align}
 \int_{\overline{\Omega}}\varphi\ (d\mu(t))(x)\geq0\qquad\text{for a.a. }t>0,\nonumber
\end{align}
so that $\mu(t)$ is for a.a. $t>0$ a non-negative measure;\\
due to \cref{conv9} and \cref{conv06} (set $\varphi\equiv 1$),  $u$ and $\mu$ satisfy  \cref{massu0}. 

Finally, passing to the limit in \cref{superve} for $k=k_m$ as $m\rightarrow\infty$, and { using \cref{convv}-\cref{convnv}, we find} that $(u,v)$ satisfies \cref{superv}. {\cref{TheoEx}} is proved.

\section{Classical solutions to \texorpdfstring{\cref{KS}}{}: proof of {\texorpdfstring{\cref{classol}}{}}}\label{SecCl}
In this final section we prove  {\cref{classol}}. 
Thus, we assume now that { the smoothness assumptions \eqref{chi1}-\eqref{uvsmooth} hold} and verify that in this case $(u,v)$ is, in fact, a classical solution to \cref{KS}. 

Let $T_0\in(0,T)$ be arbitrary but fixed. Since both $c$ and $v$ are smooth, they are uniformly bounded in $[0,{ T_0}]\times\overline{\Omega}$. From now on we assume that 
\begin{align*}
 n>\|u\|_{{ C}([0,{ T_0}]\times\overline{\Omega})}, \qquad l>\|v\|_{{ C}([0,{ T_0}]\times\overline{\Omega})}.
\end{align*}
Consequently,
\begin{align}
 \overline{u}, \overline{v}>0\qquad\text{ in }[0,{ T_0}]\times\overline{\Omega}.\label{posbar}
\end{align}

As $v$ is smooth, the weak boundary condition \cref{vweakbc} implies that $\partial_{\nu}v(t,\cdot)$ vanishes on the boundary of $\partial\Omega$ for all  $t\in(0,T_0]$. Further, since both $u$ and $v$ are smooth,  the weak boundary condition \cref{prweakbc} takes the form
\begin{align}
\left(\left(\nabla\left( \overline{u}^{a}\overline{v}^{b}\right)-\overline{u}^{a}\nabla \overline{v}^{b}\right)+\frac{1}{b}\left(b\left(M(n,l) +\overline{u}^{a}\right)-ac(u,v)\overline{u}^{a-1}\overline{v}\right)\nabla \overline{v}^{b}\right)\cdot\nu|_{\partial\Omega}=0.\nonumber
\end{align} 
This means that
\begin{align}
-a\overline{u}^{a-1}\overline{v}^{b}\left(\partial_{\nu} u-c(u,v)\partial_{\nu} v\right)-b\left(M(n,l) +\overline{u}^{a}\right)\overline{v}^{b-1}\partial_{\nu}v=0\qquad\text{ { in} }(0,{ T_0]}\times \partial\Omega.\label{prclasbc}
\end{align}
Dividing  \cref{prclasbc} by $-a\overline{u}^{a-1}\overline{v}^{b}$ and plugging the boundary condition for $v$, we conclude that $\partial_{\nu}u(t,\cdot)$ vanishes on $\partial\Omega$ for all  $t\in(0,T_0]$ as well. 

Next, exploiting the smoothness of $u$ and $v$, we integrate by parts in  \cref{superu} w.r.t. $t$ and $x$ and then apply the Du Bois-Reymond lemma. This results in the following differential inequality  in $(0,T_0]\times\Omega$:
\begin{align}
 &-a\overline{u}^{a-1}\overline{v}^{b}\left(\partial_tu -\nabla\cdot\left(\nabla u-c(u,v)\nabla v\right)\right)\nonumber\\
 &-b\left(M(n,l) +\overline{u}^{a}\right)\overline{v}^{b-1}\left(\partial_t v-(\Delta v-v+u)\right)\leq0.\label{Eq1ac}
\end{align} 
Similarly, \cref{superv} implies that $\mu=\xi\,dx$ for some  density function $0\leq\xi\in  C((0, T_0]\times\overline{\Omega})$ and
\begin{align}
 \partial_t v=\Delta v-v+u+\xi\qquad\text{in }{(0,T_0]}\times\Omega.\label{vxi}
\end{align}
Dividing \cref{Eq1ac} by $-a\overline{u}^{a-1}\overline{v}^{b}$ and making use of \cref{vxi} we deduce that
\begin{align}
 &\partial_tu -\nabla\cdot\left(\nabla u-c(u,v)\nabla v\right)\geq-\frac{b}{a}\frac{\left(M(n,l) +\overline{u}^{a}\right)}{\overline{u}^{a-1}\overline{v}}\xi\qquad\text{in }{(0,T_0]}\times\Omega.\label{Eq1ineq}
 \end{align}
 Integrating \cref{Eq1ineq} by parts over $\Omega$  using the boundary conditions and, subsequently, integrating over $(0,t)$ for any  $t\in(0,T_0]$ we then have that
\begin{align}
 \|u(t,\cdot)\|_{L^1(\Omega)}\geq& \|u_0\|_{L^1(\Omega)}-\frac{b}{a}\int_0^t\int_{\Omega}\frac{\left(M(n,l) +\overline{u}^{a}\right)}{\overline{u}^{a-1}\overline{v}}\xi\,dxds.\label{mu1}
\end{align}
On the other hand,  \cref{massu0} implies that  
\begin{align}
 \|u(t,\cdot)\|_{L^1(\Omega)}+\|\xi(t,\cdot)\|_{L^1(\Omega)}=\|u_0\|_{L^1(\Omega)}\qquad\text{ for all }{ t\in[0,T_0]}.\label{mu2}
\end{align}
Thus, combining \cref{posbar} and  \cref{mu1}-\cref{mu2}  we conclude that
\begin{align}
 \|\xi(t,\cdot)\|_{L^1(\Omega)}\leq&\frac{b}{a}\int_0^t\int_{\Omega}\frac{\left(M(n,l) +\overline{u}^{a}\right)}{\overline{u}^{a-1}\overline{v}}\xi\,dxds\nonumber\\
 \leq& \C({ T_0})\int_0^t\|\xi(s,\cdot)\|_{L^1(\Omega)}ds\qquad\text{ for all }{ t\in[0,T_0]}\nonumber
\end{align}
which yields that 
\begin{align}
\xi\equiv0\label{xi0} 
\end{align}
due to the Gronwall lemma. Plugging \cref{xi0} into \cref{vxi},  \cref{Eq1ineq}, and \cref{mu2} immediately yields that \cref{Eqv} is satisfied in the classical sense, and it holds that
\begin{align}
 &\partial_tu -\nabla\cdot\left(\nabla u-c(u,v)\nabla v\right)\geq0\qquad\text{in }{(0,T_0]}\times\Omega\label{Eq1ineq0}
 \end{align}
 and
 \begin{align}
  \|u(t,\cdot)\|_{L^1(\Omega)}=& \|u_0\|_{L^1(\Omega)}\qquad\text{for all }{ t\in[0,T_0]}.\label{mas}
 \end{align}
Since \cref{Eq1ineq0} is subject to the no-flux boundary conditions,   \cref{Eq1ineq0}  and \cref{mas} imply that equality holds in \cref{Eq1ineq0}. Thus, equation \cref{Equ} is satisfied in the classical sense, and {\cref{classol}} is proved.

\bibliographystyle{siamplain}
\bibliography{bibl}
\end{document}